\documentclass{amsart}
\begin{document}

% ----------------------------------------------------------------
\vfuzz2pt % Don't report over-full v-boxes if over-edge is small
\hfuzz2pt % Don't report over-full h-boxes if over-edge is small
%Theorems
\newtheorem{thm}{Theorem}[section]
\newtheorem{corollary}[thm]{Corollary}
\newtheorem{lemma}[thm]{Lemma}
\newtheorem{proposition}[thm]{Proposition}
\newtheorem{defn}[thm]{Definition}
\newtheorem{remark}[thm]{Remark}
\newtheorem{example}[thm]{Example}
\newtheorem{fact}[thm]{Fact}
%\numberwithin{equation}{section}
\
% MATH -----------------------------------------------------------
\newcommand{\norm}[1]{\left\Vert#1\right\Vert}
\newcommand{\abs}[1]{\left\vert#1\right\vert}
\newcommand{\set}[1]{\left\{#1\right\}}
\newcommand{\Real}{\mathbb R}
\newcommand{\eps}{\varepsilon}
\newcommand{\To}{\longrightarrow}
\newcommand{\BX}{\mathbf{B}(X)}
\newcommand{\A}{\mathcal{A}}
\newcommand{\onabla}{\overline{\nabla}}
\newcommand{\hnabla}{\hat{\nabla}}

% ----------------------------------------------------------------

\def\proof{\medskip Proof.\ }
\font\lasek=lasy10 \chardef\kwadrat="32 %kwadrat
\def\kwadracik{{\lasek\kwadrat}}
\def\koniec{\hfill\lower 2pt\hbox{\kwadracik}\medskip}

\newcommand*{\C}{\mathbf{C}}
\newcommand*{\R}{\mathbf{R}}
\newcommand*{\Z}{\mathbf {Z}}

\def\sb{f:M\longrightarrow \C ^n}
\def\det{\hbox{\rm det}\, }
\def\detc{\hbox{\rm det }_{\C}}
\def\i{\hbox{\rm i}}
\def\tr{\hbox{\rm tr}\, }
\def\rk{\hbox{\rm rk}\,}
\def\vol{\hbox{\rm vol}\,}
\def\Im {\hbox{\rm Im}\, }
\def\Re{\hbox{\rm Re}\, }
\def\interior{\hbox{\rm int}\, }
\def\e{\hbox{\rm e}}
\def\pai{\partial _i}
\def\paj{\partial _j}
\def\pu{\partial _u}
\def\pv{\partial _v}
\def\pui{\partial _{u_i}}
\def\puj{\partial _{u_j}}
\def\puk{\partial {u_k}}
\def\div{\hbox{\rm div}\,}
\def\Ric{\hbox{\rm Ric}\,}
\def\r#1{(\ref{#1})}
\def\ker{\hbox{\rm ker}\,}
\def\im{\hbox{\rm im}\, }
\def\I{\hbox{\rm I}\,}
\def\id{\hbox{\rm id}\,}
\def\exp{\hbox{{\rm exp}^{\tilde\nabla}}\.}
\def\cka{{\mathcal C}^{k,a}}
\def\ckplusja{{\mathcal C}^{k+1,a}}
\def\cja{{\mathcal C}^{1,a}}
\def\cda{{\mathcal C}^{2,a}}
\def\cta{{\mathcal C}^{3,a}}
\def\c0a{{\mathcal C}^{0,a}}
\def\f0{{\mathcal F}^{0}}
\def\fnj{{\mathcal F}^{n-1}}
\def\fn{{\mathcal F}^{n}}
\def\fnd{{\mathcal F}^{n-2}}
\def\Hn{{\mathcal H}^n}
\def\Hnj{{\mathcal H}^{n-1}}
\def\emb{\mathcal C^{\infty}_{emb}(M,N)}
\def\M{\mathcal M}
\def\Ef{\mathcal E _f}
\def\Eg{\mathcal E _g}
\def\Nf{\mathcal N _f}
\def\Ng{\mathcal N _g}
\def\Tf{\mathcal T _f}
\def\Tg{\mathcal T _g}
\def\diff{{\mathcal Diff}^{\infty}(M)}
\def\embM{\mathcal C^{\infty}_{emb}(M,M)}
\def\U1f{{\mathcal U}^1 _f}
\def\Uf{{\mathcal U} _f}
\def\Ug{{\mathcal U} _g}
\def\hnu{\hat\nu}
\def\gnu{\nu_g}
\def\C{{\mathcal ()}}
\def\A{{\mathcal (*)}}
\def\T{{\rm T}}
\title{Cauchy-Kowalevski's theorem applied for counting geometric structures}
%\thanks
\author{Barbara Opozda, W{\l}odzimierz M. Mikulski}

% Secondary: 35F35, 35G50}


\thanks{The first author   was supported by the NCN grant UMO-2013/11/B/ST1/02889}

\address{Faculty of Mathematics and Computer Science UJ, ul. {\L}ojasiewicza 6, 30-348, Cracow, Poland}

\email{Wlodzimierz.Mikulski@im.uj.edu.pl}

 \email{Barbara.Opozda@im.uj.edu.pl}

\maketitle
\begin{abstract}
How many are linear connections  with prescribed Ricci tensor? How
many are statistical structures? The questions are answered in the
analytic case  by using the Cauchy-Kowalevski theorem.
\end{abstract}

\section{Introduction}
Our study is inspired by  the recent  paper of Z. Du\v sek and O.
Kowalski \cite{DK}. Roughly speaking, the question is how many
structures of a prescribed type there exist. By a satisfactory
answer we mean a theorem saying that the set of such structures is
parametrized by some family (finite) of arbitrarily chosen
functions. We consider the local setting of the question. It turns
out that the theorem of Cauchy-Kowalevski can be used as a tool in
answering it. Of course, using this tool implies that we must
 restrict to analytic structures. But the advantage is that the tool
belongs to the fundamentals of mathematics  and a procedure of
getting structures is explicit modulo solving a  Cauchy-Kowalevski
system of differential equations. On the other hand it seems that
the method fits only  very special situations.

The paper  deals with two questions. The first one is how many
connections have a  prescribed Ricci tensor. The question of
existence of connections with prescribed Ricci tensor was studied,
for instance,  in \cite{DeT}, \cite{G1} and \cite{G}. In particular,
it was proved in \cite{G1} that any analytic  symmetric tensor of
type $(0,2)$ can be locally realized as the symmetric part of the
Ricci tensor of some torsion-free connection. We extend this result
to not necessarily symmetric prescribed tensors and the whole Ricci
tensors. Namely, we observe that a necessary condition for a tensor
of type $(0,2)$ to be (locally) the Ricci tensor of some
torsion-free connection is that its anti-symmetric part is a closed
form. For an analytic tensor field the closedness of the
anti-symmetric part is also a sufficient condition for a local
realization as the Ricci tensor of a torsion-free connection.
Moreover, we show that the set of all germs at a point in $\R^n$ of
analytic torsion-free connections $\nabla$ with prescribed Ricci
tensor (whose anti-symmetric part is closed) depends bijectively on
$\frac{n^3-3n}{2}+1$ functions of $n$ variables and
$\frac{n^2+n}{2}$  functions of $(n-1)$ variables. In particular,
the functions of $n$ variables are some Christoffel symbols of
$\nabla$. Choosing them in special ways one can produce structures
with additional properties. In the case of connections with
arbitrary torsion we prove, modifying slightly the proof of the main
theorem from \cite{DK}, that the set of all germs of connections
with prescribed Ricci tensor depends on $n^3-n^2$ analytic functions
of $n$ variables and $n^2$ functions on $n-1$ variables. We also
consider the case where the trace  of the torsion vanishes. We give
a partial result to the question how many metric structures there
are with prescribed Ricci tensor. Namely, an answer is provided in
the 2-dimensional case for non-degenerate Ricci tensors. Here  the
Cauchy-Kowalevski theorem of the second order is used.

 Another question which can be treated by means of
Cauchy-Kowalevski's theorem is the one about the amount of
statistical structures. A statistical structure is a pair
$(g,\nabla)$, where $\nabla$ is a torsion-free connection, $g$ is a
metric tensor field and $\nabla g$ as a $(0,3)$-tensor field is
symmetric. Statistical structure are examples of Codazzi pairs. Such
structures are very important in differential geometry. For
instance, the theory of equiaffine  hypersurfaces in $\R^n$ is based
on such structures. The theory of the second fundamental form of
hypersurfaces in space forms serves as another example. The induced
structures of Lagrangian submanifolds in complex space forms are
statistical structures. Statistical structures appear in statistics
and  information geometry. As regards the question we are concerned
with, we  find how many analytic statistical structures there are
around a point in $\R^n$. The proof of Theorem \ref{statistical}
provides an explicit procedure of finding  such structures.

\section{Preliminaries}

 Recall the
theorem of Cauchy-Kowalevski in the version we need for our
considerations. We adopt the notation $(f)_i=\frac{\partial
f}{\partial x^j}$, $(f)_{jk}=\frac{\partial ^2f}{\partial
x^j\partial x^k}$ for a function on a domain endowed with a
coordinate system $(x^1,...,x^n)$. All coordinate systems used in
this paper are analytic.

\begin{thm}
 Consider a system of differential equations for unknown functions
$U^1,....,U^N$ in a neighborhood of $0\in \R^n$ and of the form

\begin{eqnarray*}
&&(U^1)_1=H^1(x^1,...,x^n,U^1,...,U^N,(U^1)_2,...,(U^1)_n,...,(U^N)_2,...,(U^N)_n),
\\ &&( U^2)_1=H^2(x^1,...,x^n,U^1,...,U^N,(U^1)_2,...,(U^1)_n,...,(U^N)_2,...,(U^N)_n),
\\ && \ \ ...
\\ && (U^N)_1=H^N(x^1,...,x^n,U^1,...,U^N,(U^1)_2,...,(U^1)_n,...,(U^N)_2,...,(U^N)_n),
\end{eqnarray*}
where $H^i$, $i=1,...,N$, are  analytic functions of all variables
in a neighborhood of $(0,..., 0,\varphi ^1(0),...,\varphi^N(0),
(\varphi^1)_2(0),...,(\varphi ^1)_n(0),..., (\varphi
^N)_2(0),...,(\varphi ^N)_n(0))\in \R^{(N+1)n}$ for  analytic
functions $\varphi ^1,..., \varphi ^N$ given in a neighborhood of
$0\in \R ^{n-1}$.

Then the system has a unique solution
$(U^1(x^1,...,x^n),...,U^N(x^1,...,x^n))$ which is analytic around
$0\in \R^n$ and satisfies the initial conditions
$$
U^i(0,x^2,...,x^n)=\varphi^i(x^2,...,x^n)\ \  for \ \ i=1,...,N.
$$
\end{thm}
\bigskip

In the second order Cauchy-Kowalewski theorem we additionally
prescribe analytic functions $\psi ^1,..., \psi ^N$ defined in a
neighborhood of $0\in \R ^{n-1}$. We have $(U^1)_{11},...,
(U^N)_{11}$ on the left-hand sides and we add to the set of
arguments of $H^1,..., H^N$ on the right-hand sides the first
derivatives $(U^1)_1,...,(U^N)_1$ and the second derivatives
$(U^i)_{jk}$ for $i=1,...,N$, $j=1,...,n$ and $k=2,...,n$. To the
initial conditions we add the conditions
$$(U^i)_1(0,x^2,...,x^n)=\psi ^i(x^2,...,x^n)$$ for the prescribed
 functions $\psi^i$, $i=1,...,N$.

Since the problems we study are of local nature, we shall locate
geometric structures in open neighborhoods of $0\in\R^n$. For the
beginning a neighborhood  can be equipped with any analytic
coordinate system, for instance,  the canonical one.

 In the
following  theorems, when we write about objects in a neighborhood
of $0\in\R^n$, for instance connections, tensor fields, functions,
we mean, in fact,  their germs at $0$.

\medskip

\section{How many are connections with prescribed Ricci tensor}

For  a fixed  coordinate system $(x^1,...,x^n)$  the Ricci tensor
$\Ric$ of a linear connection $\nabla$ with  Christoffel symbols
$\Gamma^i_{jk}$ is expressed by the formula
\begin{equation}
\Ric(\partial_i,\partial_j)=\sum_{k=1}^n[(\Gamma^k_{ij})_k-
(\Gamma^k_{kj})_i]+\sum_{k,l=1}^n[\Gamma^l_{ij}\Gamma^k_{kl}-
\Gamma^l_{kj}\Gamma^k_{il}].
\end{equation}
Let $r$ be an analytic tensor field of type $(0,2)$ around
$0\in\R^n$. Set $r_{ij}=r(\partial_i,\partial_j)$.  Modifying
arguments from \cite{DK} we will prove how many real analytic linear
connections $\nabla$ exist such that $\Ric=r$.

The condition $\Ric=r$  is equivalent to the system of equations
\begin{equation}\label{RIC}
\sum_{k=1}^n[(\Gamma^k_{ij})_k-(\Gamma^k_{kj})_i]=
\sum_{k,l=1}^n[\Gamma^l_{kj}\Gamma^k_{il}-\Gamma^l_{ij}\Gamma^k_{kl}]
+r_{ij}\ , \ i,j=1,...,n.
\end{equation}
Set
\begin{equation}\Lambda
_{ij}=\sum_{k,l=1}^n[\Gamma^l_{kj}\Gamma^k_{il}-\Gamma^l_{ij}\Gamma^k_{kl}]\end{equation}
and rewrite  the system (\ref{RIC}) in the form
\begin{equation}
[(\Gamma^1_{ij})_1+...+(\Gamma^n_{ij})_n]-[(\Gamma^1_{1j})_i+...+(\Gamma^n_{nj})_i]=\Lambda_{ij}+r_{ij}\
,\  i,j=1,...,n.
\end{equation}
For  $i=1$ and $j=1,...,n$, we keep each derivative
$(\Gamma^n_{nj})_1$ on the left-hand side of the corresponding
equation. We denote the sum of all remaining terms on the left-hand
side of the corresponding equation by $\Lambda'_{1j}$ and move it to
the right-hand side. For $i>1$ and $j=1,...,n$, we keep each
derivative $(\Gamma^1_{ij})_1$ on the left-hand side of the
corresponding equation. We denote the sum of all remaining terms on
the left-hand side of the corresponding equation by $\Lambda'_{ij}$
and move it to the right-hand side. Then we obtain the (equivalent)
system
\begin{equation}\label{AA}
\begin{array}{rcl}
&&(\Gamma^n_{nj})_1=-\Lambda_{1j}-r_{1j}+\Lambda'_{1j}\ , \
j=1,...,n, \\
&&(\Gamma^1_{ij})_1=\Lambda_{ij}+r_{ij}-\Lambda'_{ij}\ , \
i=2,...,n\ , \  j=1,...,n\ . \end{array} \end{equation} We see that
the first derivatives which are on the left-hand sides of this
system are not present in any terms on the right-hand sides.
\begin{thm}\label{torsja1}
Let $r$ be an analytic tensor field of type $(0,2)$ around
$0\in\R^n$. The family of real analytic linear connections $\nabla$
with the Ricci tensor $\Ric=r$ depends bijectively on $n^3-n^2$
analytic functions of $n$ variables and $n^2$ analytic functions of
$n-1$ variables.
\end{thm}

\proof We can choose  $n^3-n^2$ Christoffel symbols $\Gamma^k_{ij}$
not present on the left hand side of (\ref{AA}) as arbitrary
analytic functions.  Then  $n^2$ analytic functions of $n-1$
variables appear by solving the system (\ref{AA}) by the
Cauchy-Kowalevski theorem. \koniec

For a linear connection $\nabla$ with  torsion $ \T
(X,Y)=\nabla_XY-\nabla_YX-[X,Y]\ , $ we have the $1$-form $\tau$
given by \begin{equation}\tau(Y)=\tr(X\to \T(X,Y))\ .
\end{equation}
Using a similar method as above, given an analytic tensor field $r$
around $0\in \R^n$, we describe all real analytic linear connections
$\Gamma$ such that $\tau=0$ and $\Ric=r$.

Clearly, this problem is equivalent to finding all solutions of the
system consisting of the system (\ref{AA}) and
\begin{equation}\label{AA1}
\sum_{i=1}^n(\Gamma^i_{ik}-\Gamma^i_{ki})=0\ , \ k=1,...,n\ .
\end{equation}

\begin{thm}\label{torsion2}
Let $n\geq 3$ and $r$ be an analytic tensor field of type $(0,2)$
around $0\in\R^n$. The family of all real analytic linear
connections $\nabla$ with $\tau=0$ and $\Ric=r$ depends bijectively
on $n^3-n^2-n$ analytic functions of $n$ variables and $n^2$
analytic functions of $n-1$ variables.
\end{thm}

\proof  From  (\ref{AA1}) we have
\begin{equation} \label{AA2}
\begin{array}{rcl}
&&\Gamma^{k+1}_{k,k+1}= -\sum_{i=1}^{k-1}\Gamma^i_{ki}-
\sum_{i=k+2}^n\Gamma^i_{ki} \\
&& \ \ \ \ +\sum_{i=1}^{k-1}\Gamma^i_{ik}
+\sum_{i=k+1}^n\Gamma^i_{ik}, \ \ \ \  \   k=1,...,n-1, \\
&& \Gamma^{n-1}_{n,n-1}=
-\sum_{i=1}^{n-2}\Gamma^i_{ni}+\sum_{i=1}^{n-1}\Gamma^i_{in}.
\end{array}
\end{equation}
 Since $n\geq 3$, the Christoffel symbols on the left-hand sides of
(\ref{AA2})  are not present on the left-hand sides of the $n^2$
equalities of (\ref{AA}). We substitute the above $n$ equalities
(\ref{AA2}) into the  $n^2$ equalities of (\ref{AA}). We obtain
\begin{equation}\label{AA3}
\begin{array}{rcl}
&& (\Gamma^n_{nj})_1=-\tilde\Lambda_{1j}-r_{ij}+\tilde\Lambda'_{1j}\
,\ j=1,...,n\ ,\\
&&(\Gamma^1_{ij})_1=\tilde\Lambda_{ij}+r_{ij}-\tilde\Lambda'_{ij}\ ,
i=2,...,n\ , j=1,...,n\ ,
\end{array}
\end{equation}
where  $\tilde\Lambda_{1j}$, $\tilde\Lambda'_{1j}$,
$\tilde\Lambda_{ij}$, $\tilde\Lambda'_{ij}$ are $\Lambda_{1j}$,
$\Lambda'_{1j}$, $\Lambda_{ij}$, $\Lambda'_{ij}$ respectively, after
the substitutions. It is easy to see that  the first derivatives
which are on the left-hand sides of the system (\ref{AA3}) are not
present on the right-hand sides. Now we can choose $n^3-n^2-n$
Christoffel symbols $\Gamma^{i}_{jk}$ not present on the left hand
sides of (\ref{AA3}) and of (\ref{AA2}) as arbitrary analytic
functions. Then $n^2$ analytic functions of $n-1$ variables appear
by solving (\ref{AA3}) by means of the Cauchy-Kowalevski theorem.
\koniec

\medskip

 If $n=2$ then the condition  $\tau=0$ yields $\T=0$. Hence the connection is torsion-free. We shall now
 study this case.
Set
\begin{equation}
D_j=div ^\nabla\paj =\tr(X\to\nabla _X\paj)=\sum_{k=1}^n\Gamma
^k_{kj}.
\end{equation}
Then the formula for the Ricci tensor can be written as  follows

\begin{equation}\Ric (\pai,\paj)=\sum_{k=1}^n(\Gamma
_{ij}^k)_k -(D_j)_i +\Lambda _{ij}.\end{equation}

We decompose the Ricci tensor into its symmetric and anti-symmetric
parts, that is, $\Ric =s+a$, where
\begin{equation}
s(X,Y)=\frac{\Ric (X,Y) +\Ric(Y,X) }{2},\ \ \ \ \ a(X,Y)=\frac{\Ric
(X,Y)-\Ric(Y,X)}{2}.
\end{equation}
For a torsion-free connections the portions $\sum_{k=1}^n(\Gamma
^k_{ij})_k$ and $\Lambda _{ij}$ are symmetric for $i$ and $j$. Hence
for a torsion-free connection we have
\begin{equation}\label{a}
a_{ij}=a(\pai,\paj)=\frac{(D_i)_j-(D_j)_i}{2},
\end{equation}
\begin{equation}\label{s}
s_{ij}= s(\pai,\paj)=\sum_{k=1}^n(\Gamma^k_{ij})_k
-\frac{(D_j)_i+(D_i)_j}{2}+\Lambda_{ij}.
\end{equation}

In \cite{O}   the following proposition was proved. Since its proof
is short, we cite it here.
\begin{proposition}
For a torsion-free connection on a paracompact manifold $M$ the
anti-symmetric part of its Ricci tensor  is exact.
\end{proposition}
\proof By the first Bianchi identity we have
$$
\tr R(X,Y)=\Ric(Y,X)-\Ric (X,Y)
$$
for a torsion-free connection $\nabla$, where $R$ is its curvature
tensor. Let $\nabla'$ be any torsion-free connection whose Ricci
tensor $\Ric '$ is symmetric. It can be the Levi-Civita connection
of some metric. Denote by $Q$  the difference tensor between
$\nabla$ and $\nabla '$, that is, $Q(X,Y)=Q_XY=\nabla
_XY-\nabla'_XY$. Define the 1-form $\delta$ on $M$ by
$$
\delta(X)=\tr Q_X.
$$
Then
$$
d\delta (X,Y)=\frac{1}{2}\{\tr \nabla' Q(X,Y,\cdot)-\tr \nabla'
Q(Y,X,\cdot)\}.
$$
The curvature tensors $R$ and $R'$ for $\nabla$ and $\nabla '$ are
related by the formula
$$
R(X,Y)Z=R'(X,Y)Z +\nabla' Q(X,Y,Z)-\nabla
'Q(Y,X,Z)+Q_XQ_YZ-Q_YQ_XZ.$$ It follows that $\tr R(X,Y)=\tr
R'(X,Y)+2d\delta (X,Y)=2d\delta (X,Y)$. \koniec

Since we study problems of local nature, we replace the exactness of
the form in the above theorem by its closedness. We shall prove
\begin{thm}\label{torsion-free}
A real analytic tensor field $r$ of type $(0,2)$ can be locally
realized as the Ricci tensor of a torsion-free connection if and
only if its anti-symmetric part $a$, that is,
$a(X,Y)=\frac{r(X,Y)-r(Y,X)}{2}$, is closed. For a given tensor
field $r$ in a neighborhood of  $0\in \R^n$ satisfying the above
conditions the set of    all analytic torsion-free connections whose
Ricci tensor is $r$, depends bijectively on $\frac{n^3-3n}{2}+1$
arbitrarily chosen analytic functions of $n$ variables and
$\frac{n^2+n}{2}$ arbitrarily chosen analytic functions of $n-1$
variables.\end{thm}

\proof Let $s$ denote the symmetric part of $r$. The functions
$a_{ij}=a(\pai,\paj)$, $s_{ij}=s(\pai,\paj)$ are given. Assume that
the form $a$ is closed and $r$ is analytic. It is locally exact,
hence around the fixed point $0$  there is an analytic $1$-form
$\alpha$ such that $a=-d\alpha$. The $1$-form $\alpha$ is chosen up
to one function, that is, $\alpha$ can be replaced by $\alpha
+d\phi$ for any function $\phi$. Let $\alpha =\alpha_1dx^1 +...
+\alpha _n dx^n$.
%The functions $\alpha_1,..., \alpha_n$ are analytic.
We have $2a_{ij} =-2d\alpha(\pai,\paj)=(\alpha_i)_j -(\alpha_j)_i$.
Suppose that $r$ is the Ricci tensor of some torsion-free connection
whose Christoffel symbols $\Gamma _{ij}^k$ are unknown.  Then
\begin{equation}\label{aD}
\frac{(D_i)_j +(D_j)_i}{2} =a_{ij}+(D_j)_i\end{equation} for
$i,j=1,...,n$. Set $D_i=\alpha_i $ for $i=1,,,.,n$. We have already
used (\ref{a}) and from now on the functions $D_1,...,D_n$ are
given.

All the conditions from (\ref{s}) must be satisfied. We have
$$s_{11}=\sum _{k=1}^n (\Gamma^k_{11})_k -(D_1)_1+\Lambda_{11},$$
hence
$$(\Gamma_{11}^1+\Gamma^2_{21}+...+\Gamma^n_{n1})_1=(\Gamma^{1}_{11})_1+(\Gamma
^2_{11})_2 +...+(\Gamma ^n_{11})_n +\Lambda _{11}-s_{11}.$$ We can
write it equivalently as
\begin{equation}\label{1}
(\Gamma^2_{12})_1=\sum_{k=2}^n(\Gamma^k_{11})_k
-\sum_{k=3}^n(\Gamma^k_{k1})_1 +\Lambda_{11}-r_{11}.\end{equation}
For $i>1$ we have
$$
s_{1i}=\sum_{k=1}^n(\Gamma
^k_{1i})_k-\frac{(D_i)_1+(D_1)_i}{2}+\Lambda_{1i}.
$$
By using (\ref{aD}) we get
$$(\Gamma^1_{1i})_1=-(\Gamma^2_{1i})
_2 -...-(\Gamma ^n_{1i})_n -\Lambda _{1i}+a_{i1}+(D_1)_i+s_{i1}.$$
We  can write it as follows
\begin{equation}\label{2}
(\Gamma^1_{1i})_1=-(\Gamma^2_{1i})_2-...-(\Gamma^n_{1i})_n-\Lambda_{1i}+(D_1)_i+r_{i1}.
\end{equation}
For $i,j$, where $1<i\le j\le n$, we have
$$s_{ij}=\sum_{k=1}^n(\Gamma^k_{ij})_k
-\frac{(D_j)_i+(D_i)_j}{2}+\Lambda _{ij},$$ that is,
$$s_{ij}=(\Gamma^1_{ij})_1+(\Gamma^2_{ij})_2+...+(\Gamma
^n_{ij})_n-a_{ij}-(D_j)_i +\Lambda_{ij}.$$ We shall write it as
follows
\begin{equation}\label{3}
(\Gamma ^1_{ij})_1=-(\Gamma ^2_{ij})_2-...-(\Gamma^n_{ij})_n-\Lambda
_{ij}+(D_j)_i +r_{ij}.
\end{equation}
Collecting the equations from (\ref{1})-(\ref{3}) we get the
following Cauchy-Kowalevski system of $\frac{n(n+1)}{2}$ equations
(equivalent to (\ref{s}))

\begin{equation}\label{CK}
\begin{array}{rcl}
&&(\Gamma^2_{12})_1=\sum_{k=2}^n(\Gamma^k_{11})_k
-\sum_{k=3}^n(\Gamma^k_{1k})_1 +\Lambda_{11}-r_{11},\\
&&
(\Gamma^1_{1i})_1=-(\Gamma^2_{1i})_2-...-(\Gamma^n_{1i})_n-\Lambda_{1i}+(D_1)_i+r_{i1},\ \ \ \ \ i>1,\\
&&(\Gamma ^1_{ij})_1=-(\Gamma
^2_{ij})_2-...-(\Gamma^n_{ij})_n-\Lambda _{ij}+(D_j)_i +r_{ij},  \ \
\ \ \  1<i\le j\le n.
\end{array}
\end{equation}
The quantities  $r_{11}$, $(D_1)_i +r_{i1}$, $(D_j)_i +r_{ij}$ are
given.

Except for the dependence given by (\ref{CK}) the Christoffel
symbols are related by the following  system of equations
\begin{equation}\label{D}
\begin{array}{rcl}
&&D_1=\Gamma^1_{11}+[\Gamma^2_{21}]+...+\Gamma^n_{n1},\\
&&D_2=[\Gamma^1_{12}]+\Gamma^2_{22}+...+\Gamma^n_{n2},\\
&&\cdot\\
&&\cdot\\
&&\cdot\\
&&D_n=[\Gamma^1_{1n}]+\Gamma^2_{2n}+...+\Gamma^n_{nn},
\end{array}
\end{equation}
where, by using brackets, we marked
 the  Christoffel symbols from the right-hand side of (\ref{D}) which
 appear on the left-hand side of (\ref{CK}). Observe also that on
 the right-hand sides of (\ref{D}) there are no Christoffel symbols
 which repeat because of the symmetry of $\Gamma^k_{ij}$ in lower
 indices.

 From each of the equations in (\ref{D}) we want to determine one
 Christoffel symbol and then substitute it into   (\ref{CK}) by the expression obtained from (\ref{D}). Of
 course, we should not  determine and substitute any marked symbol. Moreover, we have to
 do it in such a way that, after the substitution into  (\ref{CK}),
 the derivatives from the left-hand side of  (\ref{CK}) will not
 appear on the right-hand side of  (\ref{CK}). Therefore, from the first equation
 of (\ref{D}) we can only take:
 $\Gamma^1_{11}=D_1-\Gamma^2_{21}-...-\Gamma^n_{n1}$.
 From the next equations we  can take $\Gamma ^k_{kk}$ (but here it is
 not necessary to do it in this way).

 For the modified system (\ref{CK}) (after  the substitutions) we
 can apply the Cauchy-Kowalevski theorem.

 \medskip

 We shall now count how many Christoffel symbols can be chosen
 arbitrarily.
Note that all Christoffel symbols for which the upper index is equal
to one or two of lower indices are on the right-hand side of
(\ref{D}). We see that from (\ref{D}) we can choose $n(n-2)$ symbols
arbitrarily. Consider now the Christoffel symbols for which the
upper index is different than each of the lower indices. Consider
first the symbols whose upper index is 1. All of them appear on the
left-hand side of (\ref{CK}) so we cannot choose them. Finally
consider those Christoffel symbols whose upper index is $k$, where
$1<k\le n$, and $k$ is different than any of the lower indices. They
do not appear neither on the left-hand side of (\ref{CK}) nor on the
right-hand side of (\ref{D}). All of them can be chosen arbitrarily.
There are $(n-1)\frac{(n-1)n}{2}=\frac{(n-1)^2n}{2}$ such symbols.
Therefore we can choose $n(n-2)+\frac{(n-1)^2n}{2}=\frac{n^3-3n}{2}$
Christoffel symbols arbitrarily. The function $\phi$ from the
beginning of the proof is also an arbitrarily chosen function of $n$
variables.

\koniec

\begin{remark}
{\rm For $n=2$ we have
\begin{equation}\label
{n=2}
\begin{array}{rcl}
&&D_1=\Gamma^1_{11}+[\Gamma^2_{12}],\\
&&D_2=[\Gamma^1_{21}]+\Gamma^2_{22}.\\
\end{array}
\end{equation}
None of the Christoffel symbols from the right-hand side of
(\ref{n=2}) can be chosen arbitrarily (in the above procedure). We
have $\frac{n^3-3n}{2}=1$. The only Christoffel symbol  which can be
 arbitrarily chosen in this case is $\Gamma^2_{11}$. In particular, we can choose it 0 and then
the vector field $\nabla _{\partial_1}\partial _1$ is parallel to
$\partial _1$ (but  we cannot assume that $\nabla
_{\partial_1}\partial _1$ vanishes). For any dimension   the
functions $\Gamma^k_{11}$ for $k=2,...,n$ are up to choice. In
particular, one can choose them 0, which means that $\nabla
_{\partial_1}\partial _1$ is parallel to $\partial _1$. But we
cannot assume that $\Gamma ^1_{11}=0$. From the last equation of
(\ref{CK}) it is clear that we cannot assume that for some $i>1$ we
have $\Gamma ^k_{ii}=0$ for all indices $k$, because we cannot
choose $\Gamma ^1_{ii}$ arbitrarily.}

\end{remark}

\bigskip

We shall now  give a partial answer to the question how many
Levi-Civita connections are those whose Ricci tensor is a prescribed
symmetric tensor $r$ of type $(0,2)$.

For a metric tensor field $g$ (not necessarily positive definite)
   the Christoffel symbols of its Levi-Civita connection are
given by $$
\Gamma^s_{ij}=\frac{1}{2}\sum_{k=1}^ng^{sk}\left((g_{ki})_j
+(g_{jk})_i-(g_{ji})_k\right).$$ where
$g_{ij}=g(\partial_i,\partial_j)$ and $(g^{sk})$ is the inverse
matrix of the matrix $(g_{ij})$.
 If $n=2$ and the matrix $(g_{ij})$  has a diagonal form in the coordinate system then
\begin{equation}\label{Gamma_dla_g}
\begin{array}{crl}
&&\Gamma^1_{11}=\frac{1}{2}g^{11}( g_{11})_1, \
\Gamma^2_{11}=-\frac{1}{2}g^{22}( g_{11})_2,\
\Gamma^1_{21}=\Gamma^1_{12}=\frac{1}{2}g^{11}( g_{11})_2,\\
&& \Gamma^2_{21}=\Gamma^2_{12}=\frac{1}{2}g^{22}( g_{22})_1, \
\Gamma^1_{22}=-\frac{1}{2}g^{11}( g_{22})_1,\ \
\Gamma^2_{22}=\frac{1}{2}g^{22}( g_{22})_2,
\end{array}
\end{equation}
where $g^{11}={1\over g_{11}}$ and $g^{22}={1\over g_{22}}$. The
Ricci tensor $\Ric$ of the Levi-Civita connection for $g$ satisfies
the equality $\Ric =fg$, where $f$ is the sectional curvature of
$g$. Using  (\ref{Gamma_dla_g}), by a straightforward computation
one gets
\begin{equation}\label{sectional_curvature}
\begin{array} {rcl}
&&f= -\frac{1}{2}g^{11}g^{22}[(g_{11})_{22}+(g_{22})_{11}]\\&&\ \ \
\ \  \ \  +\frac{1}{4}
g^{11}(g^{22})^2[(g_{22})_2(g_{11})_2+((g_{22})_1)^2]\\&&\ \ \ \ \ \
\  + \frac{1}{4} (g^{11})^2g^{22}[( g_{11})_1(g_{22})_1+(
(g_{11})_2)^2].
\end{array}
\end{equation}

Note that for an analytic metric tensor field on a 2-dimensional
manifold there is an analytic orthogonal coordinate system around
each point of the domain of the metric tensor field.

\begin{thm}\label{metric}
Let $r$ be an analytic non-degenerate  tensor field of type $(0,2)$
such that its matrix is diagonal in an analytic coordinate system
$(x_1,x^2)$ on a neighborhood of $0\in\R^2$. Then the set of all
analytic metric tensor fields such that their Ricci tensors equal to
$r$ depends bijectively on arbitrarily chosen pairs $(\varphi,\psi)$
of analytic functions of one variable with $\varphi (0)\not=0$.
\end{thm}

\proof Suppose that $g$ is an analytic  metric  tensor field around
$0\in\R^2$ such that its Ricci tensor  $\Ric$ is equal to $r$. Then
$g=hr$ for some analytic map $h$  around $0\in\R^2$ with
$h(0)\not=0$. By (\ref{sectional_curvature}) the equality $\Ric=r$
is equivalent to the partial differential equation
\begin{equation}
\begin{array}{rcl}
&&-\frac{1}{2hr_{22}}[(hr_{11})_{22}+ (h r_{22})_{11}]\\
&&+\frac{1}{4(hr_{22})^2}[(hr_{22})_2 (hr_{11})_2+((hr_{22})_1)^2]\\
&& +\frac{1}{4h^2r_{11}r_{22}}[ (hr_{11})_1(hr_{22})_1+((h
r_{11})_2)^2]=r_{11}.
\end{array}
\end{equation}
Applying the Leibniz rule one sees that this equation can be
transformed equivalently into the one of the form
$$(h)_{11}=F(h, (h)_1, (h)_2,(h)_{12}, (h)_{22})$$ for  some
analytic map $F$.  Our theorem now follows from the
Cauchy-Kowalevski theorem of  order 2, where  two analytic functions
$\varphi$, $\psi$ of one variable are prescribed and the initial
conditions are: $h(0,x^2)=\varphi$,  $(h)_1(0,x^2)=\psi$. \koniec

\bigskip
\section{How many are statistical structures}

\bigskip
Recall that a statistical structure  on a manifold $M$ is a pair
$(g,\nabla)$, where $g$ is a  metric tensor field and $\nabla$ is a
torsion-free connection on $M$ satisfying the Coddazzi condition
saying that $\nabla g$ as a  cubic form is totally symmetric. We
assume that the metric is positive definite. A statistical structure
is called trace-free when the volume form $\nu _g$ determined by $g$
is parallel relative to $\nabla$. The trace-free statistical
structures correspond to Blaschke structures in affine differential
geometry and to minimal submanifolds in the theory of Lagrangian
submanifolds. We begin with the 2-dimensional case.

In the following Theorems \ref{statistical2} and \ref{statistical}
the metric tensor fields are unknowns, but according to the
Cauchy-Kowalevski theorem they can be arbitrarily chosen at the
point $0$. Up to linear isomorphism of $\R ^n$ we can assume that
the matrix $(g_{ij})$ at $0$ is the identity one. We  make this
assumption for both Theorems \ref{statistical2},  \ref{statistical},
that is, we assume that for the functions $g_{ij}$ appearing in
these theorems $g_{ij}(0)=\delta_{ij}$.

\begin{proposition}\label{statistical2}
For  any analytic  linear connection $\nabla$ in a neighborhood of
$0\in \R^2$ there is an analytic metric tensor field $g$
% (of a prescribed algebraic type; for instance, positive definite)
 around
$0$ such that the cubic form $\nabla g$ is symmetric. The set of all
such metric tensor fields depends on one function $g_{11}$ of two
variables and two functions $g_{12}$, $g_{22}$ of one variable.  If
additionally the Ricci tensor of $\nabla$ is symmetric then there is
a  metric tensor field $g$ such that $(g, \nabla) $ is a trace-free
statistical structure. The set of such metric tensor fields depends
on the two functions $g_{12}$, $g_{22}$ of one variable.
%The set of germs of these structures depends on two
%functions of one variable.
\end{proposition}
\proof For a metric tensor field  $g=g_{ij}$ the $(0,3)$-tensor
$\nabla g$ is symmetric if and only if
$$\nabla g(\partial _1,\partial
_2,\partial _1) =\nabla g(\partial _2,\partial_1,\partial _1),\ \ \
\ \  \ \nabla g(\partial _1,\partial _2,\partial _2) =\nabla
g(\partial _2,\partial_1,\partial _2).$$ It leads to the following
Cauchy-Kowalevski system  of differential equations with unknowns
$g_{12}$, $g_{22}$

\begin{equation}\label{2-wym}
\begin{array}{rcl}
&&(g_{12})_1=(g_{11})_2+\Gamma^1_{11}g_{12}+\Gamma^2_{11}g_{22}-\Gamma
^1_{21}g_{11}-\Gamma ^2_{21}g_{21},\\
&&(g_{22})_1=(g_{12})_2+\Gamma ^1_{12}g_{12}+\Gamma
^2_{12}g_{22}-\Gamma ^1_{22}g_{11}-\Gamma ^2_{22} g_{12}.
\end{array}
\end{equation}

The function $g_{11}$ can be arbitrary (modulo the assumption made
before the theorem). Assume now that the Ricci tensor of $\nabla$ is
symmetric. In a neighborhood of $0$ there is a volume form $\nu$
such that $\nabla \nu =0$. We want to have $\nu=\nu _g$ (up to some
constant c), that is,
\begin{equation}\label{nu=nu_g}
{\rm c}\nu (\partial _1,\partial _2)^2=g_{11}g_{22} -g_{12}^2.
\end{equation} Since $g_{22}\ne 0$  at $0$, we can
determine $g_{11}$ from (\ref{nu=nu_g}) and make the substitution
into (\ref{2-wym}). After the substitution the system remains
solvable with two arbitrarily prescribed (modulo the assumption made
before the theorem) functions of one variable.\koniec

The above consideration cannot be repeated  in  more dimensional
cases. But we have
\begin{thm}\label{statistical}
The set  of all analytic statistical structures $(g,\nabla)$ around
$0\in \R ^n$, where $n>2$,  depends on $\frac{n^3+6n^2+5n}{6}$
arbitrarily chosen analytic functions of $n$ variables, from which
one function is $g_{11}$ and $\frac{n^3+6n^2+5n-6}{6}$ functions are
some Christoffel symbols of $\nabla$, and $ \frac{n(n+1)}{2}-1$
arbitrarily chosen analytic
 functions $g_{ij}$, for $(ij)\ne (1,1)$, of $(n-1)$ variables.
\end{thm}

\proof We shall need the following lemma \vskip0.2cm

\begin{lemma}   A pair $(g,\nabla)$ is a statistical structure if and
only if
\begin{equation}\label{symmetry}
(\nabla g)(\partial_i,\partial_j,\partial_k)=(\nabla
g)(\partial_j,\partial_i,\partial_k)
\end{equation}
 for every $i,j,k=1,...,n$ with $i< j$ and $i\leq k$.\end{lemma}

{\bf \proof} The symmetry for the last two arguments of $\nabla g$
holds because  of the symmetry of $g$. Assume (\ref{symmetry}) for
all $i,j,k$ with $i<j$ and $i\leq k$. Take $i,j,k\in \{1,...,n\}$
such that $i<j$ and $k<i$. Hence $k<j$. We now have
\begin{eqnarray*}
&&(\nabla g)(\partial_i,\partial_j,\partial_k)=\nabla
g(\partial_i,\partial_k,\partial_j)=(\nabla
g)(\partial_k,\partial_i,\partial_j)\\
&& \ \ \ = (\nabla g)(\partial_k,\partial_j,\partial_i)=(\nabla
g)(\partial_j,\partial_k,\partial_i)=(\nabla
g)(\partial_j,\partial_i,\partial_k).
\end{eqnarray*}
\koniec

\medskip

 Consider first the conditions
(\ref{symmetry}) for the indices $1,k,j$, where $1\leq k\leq j$,
that is,
$$(\nabla g)(\partial_1,\partial_j,\partial_k)=(\nabla g)(\partial_j,\partial_1,\partial_k).$$
The conditions lead to the equations
\begin{equation}\label{CK-g}
(g_{jk})_1=(g_{1k})_j+\sum _{l=1}^ng_{jl}\Gamma ^l_{1k}-\sum_{l=1}^n
g_{1l}\Gamma ^l_{jk}.
\end{equation}
There are $\frac{n(n+1)}{2}-1$ equations in (\ref{CK-g}). The
portion $-1$ comes from $(g_{11})_1$. The system (\ref{CK-g}) will
be our Cauchy-Kowalevski system. According to the Cauchy-Kowalevski
theorem we can prescribe all functions $g_{jk}$ at $0$ (in
particular). We choose them such that the matrix $g_{jk}(0)$ is the
identity one. The function $g_{11}$ can be chosen arbitrarily modulo
the assumption that $g_{11}(0)=1$.

We now take into account the conditions
\begin{equation}
(\nabla g)(\partial _1,\partial _k,\partial _j)= (\nabla g)(\partial
_k,\partial _1,\partial _j)
\end{equation}
for $1<k<j$.  The conditions are equivalent to the equalities
\begin{equation}
(g_{kj})_1-\sum _{l=1}^ng_{kl}\Gamma ^l_{1j}=(g_{1j})_k-\sum_{l=1}^n
g_{1l}\Gamma ^l_{kj}.
\end{equation}
Since we postulate that $g_{jk}=g_{kj}$ and $\Gamma ^l_{kj}=\Gamma
^l_{jk}$, by using  (\ref{CK-g}) we get the conditions
\begin{equation}\label{gamma_k1j}
(g_{1k})_j+\sum _{l=1}^ng_{jl}\Gamma ^l_{1k} =(g_{1j})_k+\sum
_{l=1}^ng_{kl}\Gamma ^l_{1j}
\end{equation}
for $1<k<j\leq n$.
%, where we postulate $g_{ab}=g_{ba}$.
 We have
$(n-2)+(n-3)+...+1= \frac{(n-1)(n-2)}{2}$  equalities in
(\ref{gamma_k1j}). To each equality of (\ref{gamma_k1j}) we assign
the unique pair $(k,j)$
%of indices
%on $(g_{1k})_j$ or on $(g_{1j})_k$
with $k<j$. The obtained correspondence is a bijection between the
set of equalities (\ref{gamma_k1j}) and the set of pairs $(k,j)$ of
integers with $1<k<j\leq n$. So, we can order the system
(\ref{gamma_k1j}) by the inverse lexicographic ordering in pairs
$(k,j)$, that is, $(k_1,j_1)\leq (k,j)$ if and only if  $j< j_1$ or
$j=j_1$ and $k\leq k_1$.

 The rest of the conditions from (\ref{symmetry}) deal
with $(\nabla g)(\partial_i,\partial_j,\partial_k)$, where all $i,
j,k$ are different than $1$. Assume first that  two of the indices
$i,j,k$ are equal. We have the equalities
$$(\nabla g)(\partial_i,\partial_j,\partial_i)=(\nabla g)(\partial_j,\partial_i,\partial_i),$$
where $i=2,...,n$ and $j\in\{2,...,n\}\setminus\{i\}$. We have
$(n-1)(n-2)$  equalities here. They lead to the  conditions
\begin{equation}\label{condition_ij}
(g_{ji})_i-\sum _{l=1}^n g_{jl}\Gamma ^l_{ii}=(g_{ii})_j- \sum
_{l=1}^ng_{il}\Gamma ^l_{ji}
\end{equation}
for $i=2,...,n$ and $j\in\{2,...,n\}\setminus\{i\}$.
%, where we
%postulate $g_{ab}=g_{ba}$ and $\Gamma^a_{bc}=\Gamma^a_{cb}$.
 To
each equation from (\ref{condition_ij}) we assign the unique pair
$(j,i)$ of indices. The obtained correspondence is a bijection
between the set of equalities (\ref{condition_ij}) and the set of
pairs $(j,i)$ such that $i=2,...,n$ and
$j\in\{2,...,n\}\setminus\{i\}$. We order the system
(\ref{condition_ij}) by means of the inverse lexicographic ordering
in pairs $(i,j)$.

Consider now (\ref{symmetry}) for all remaining $i,j,k$, that is,
for $i,j,k$ such that $2\leq i<j\leq n$ and
$k\in\{2,...,n\}\setminus\{i,j\}$ and $i\leq k$. In fact $i<k$.
Hence the condition (\ref{symmetry}) gives here $$\sum
_{i=2}^{n-2}(n-i)(n-i-1)=\sum_{l=1}^{n-2}l(l-1)=\frac{n^3-6n^2+11n-6}{3}$$
equalities
\begin{equation}\label{condition_ijk}
(g_{jk})_i-\sum _{l=1}^n g_{jl}\Gamma ^l _{ik}=(g_{ik})_j- \sum
_{l=1}^n g_{il}\Gamma ^l_{jk}
\end{equation}
for $2\leq i<j\leq n$ and $k\in\{2,...,n\}\setminus\{i,j\}$ and
$i\leq k$. To each equality from (\ref{condition_ijk}) we assign the
unique triple $(i,j,k)$ of indices. This correspondence is a
bijection between the set of equalities in (\ref{condition_ijk}) and
the set of triples $(i,j,k)$ of integers such that $2\leq i<j\leq
n$, $k\in\{2,...,n\}\setminus\{i,j\}$ and $i\leq k$.  We can order
the equalities in (\ref{condition_ijk}) by means of the inverse
lexicographic ordering in triples $(i,j,k)$.

Denote by   $\A$ the ordered system of algebraic equations with
unknown Christoffel symbols consisting of the above ordered systems
(\ref{gamma_k1j}), (\ref{condition_ij}) and (\ref{condition_ijk}) in
the sequence (\ref{gamma_k1j}), (\ref{condition_ij}),
(\ref{condition_ijk}) . From each equation of the system  $\A$,
starting from the first equation and going up to the last equation,
we want to determine one Christoffel symbol and substitute it into
the Cauchy-Kowalevski system as well as into all next equations from
our system $\A$. At each step of the procedure the equations in our
Cauchy-Kowalevski system will change and the algebraic equations
will change as well. From the subsystem (\ref{gamma_k1j}) we shall
determine symbols $\Gamma^j_{1k}$, from the subsystem
(\ref{condition_ij}) we shall determine $\Gamma ^j_{ii}$ and from
the last subsystem (\ref{condition_ijk}) the symbols
$\Gamma^j_{ik}$. At each step of the procedure our system of
differential equations will remain a Cauchy-Kowalevski system and
the coefficient in front of the symbol which will be determined at a
consecutive step will be non-zero (in some neighborhood of the point
$0$), that is, it will be possible to determine this symbol from the
equation. To this aim we  assumed that the  matrix of $g$ at $0$ is
the identity one. Namely, the  system $\A$  evaluated at $0$ is the
following
\begin{equation}
\begin{array}{rcl}\label{at_p}
&&(g_{1k})_j+\Gamma ^j_{1k}=(g_{1j})_k+\Gamma^k_{1j}\ \
\ \ {\rm for}\  1<k<j\le n,\\
&& (g_{ji})_i-\Gamma ^j_{ii}=(g_{ii})_j-\Gamma^i_{ji}\ \
\ \  {\rm for}\  i=2,...,n; \ j\in\{2,...,n\}\setminus\{i\},\\
&& (g_{jk})_i-\Gamma ^j_{ik}=(g_{ik})_j-\Gamma^i_{jk}\ \
\ \  {\rm for}\  2\le i<j\le n; k\ne j; i< k.\\
\end{array}
\end{equation}

Each Christoffel symbol which we want to determine from the system
(\ref{at_p}) appears in the system only once. It is  easily seen
that we can safely apply the procedure described above to
(\ref{at_p}) because the coefficient in front of a Christoffel
symbol which we want to determine at a certain step is non-zero.
Hence it is non-zero around a point $0$ (because at each step of our
procedure we use only elementary algebraic operations) and we can
determine this Christoffel symbol in a neighborhood of $0$.

After solving the Cauchy-Kowalevski system one goes back to the
algebraic system and going from the last to the first equation one
gets a complete set of Christoffel symbols.

We have  presented an explicit procedure of solving  the system of
algebraic equations. But, if one does not want an explicit
procedure, one can shortly argue as follows. We have the system of
algebraic equations with unknowns being Christoffel symbols and
consisting of (\ref{gamma_k1j}) for $1<k<j\le n$,
(\ref{condition_ij}) for $i=2,...,n; \
j\in\{2,...,n\}\setminus\{i\}$ and (\ref{condition_ijk}) for $2\le
i<j\le n; k\ne j; i< k$. The system has the form (\ref{at_p}) at
$0$. The matrix of the coefficients of the system is of maximal rank
at $0$ and so it is around $0$. Hence the system has an analytic
solution around  the point $0$ depending on
 $$\frac{n^2(n+1)}{2}-\frac{(n-1)(n-2)}{2}-(n-1)(n-2)-\frac{n^3-6n^2+11n-6}{3}=\frac{n^3+6n^2+5n-6}{6}$$
arbitrarily chosen analytic parameters. It is seen that the
substitution of  the solutions into the system (\ref{CK-g}) does not
destroy its property of being a Cauchy-Kowalevski system.\koniec

\begin{remark}{\rm The above theorem and its proof are also valid for $g$ being pseudo-Riemannian
metric tensor fields with a fixed signature (modulo  linear
isomorphisms,  as in the above theorem). In the proof it is
sufficient to take the matrix of $g(0)$ in an appropriate form.}
\end{remark}


\bigskip
\bigskip
\begin{thebibliography}{20}

\bibitem{DeT} DeTurck D., \emph{Existence of metrics with prescribed Ricci curvature: Local
theory}, Invent. Math. 65, 1981, 179-207.


%\bibitem{DeT1} D. DeTurk, \emph{Matrics wih prescribed Ricci
%curvature}, Seminar on Differential Geometry, S.T.Yau ed., Annals of
%Math. Study 102, Princeton U. Press, 1982, 525-537.


\bibitem{DK}  Du\v sek Z.,  Kowalski O., \emph{How many are Ricci flat affine
connections with arbitrary torsion}, preprint, 2015.

\bibitem{G} Gasqui J., \emph{Connexions \`a courbure de Ricci donn\'ee}, Math.
Z., 168, 1975, 167-179.

\bibitem{G1} Gasqui J., \emph{Sur la courbure de Ricci d'une connexion lin\'eaire.}
C. R. Acad. Sc. Paris, 281, 1975, 389-391.


\bibitem{O} Opozda, B., \emph{On some properties of the curvature and Ricci
tensors in complex affine  geometry}, Geom. Dedic. 55, 1995,
141-163.

\end{thebibliography}
\end{document}